\documentclass[12pt,a4paper]{amsart}
\usepackage{amssymb,amsmath,amsthm}
\usepackage{epsfig}

\newtheorem{theorem}{Theorem}[section]
\newtheorem{lemma}{Lemma}[section]

\newtheorem{prop}{Proposition}[section]
\newtheorem{claim}{Claim}[section]

\newcommand{\al}{\alpha}
\newcommand{\be}{\beta}

\newcommand{\eps}{\varepsilon}
\newcommand{\R}{\mathbb{R}}

\newcommand{\conv}{\mathrm{conv}}

\newcommand{\p}{\mathcal{P}}
\newcommand{\aff}{\textrm{aff\;}}
\newcommand{\rank}{\textrm{rank\;}}

\numberwithin{equation}{section}

\begin{document}

\title{Tverberg's theorem, a new proof}
\author{Imre B\'ar\'any}

\begin{abstract} We give a new proof Tverberg's famous theorem: For every set $X \subset \R^d$ with $|X|=(r-1)(d+1)+1$, there is a partition of $X$ into $r$ sets $X_1,\ldots,X_r$ such that $\bigcap_{p=1}^r \conv X_p\ne \emptyset$. The new proof uses linear algebra, specially structured matrices, the theory of linear equations, and Tverberg's original ``moving the points" method.
\end{abstract}

\maketitle

\section{Introduction}\label{sec:introd}

A finite set $X \subset \R^d$ has, by definition, a Tverberg partition into $r$ parts $X_1,\ldots,X_r$ if $\bigcap_{p=1}^r \conv X_p$ is nonempty, of course $X_1,\ldots,X_r$ form a partition of $X$. In this case $X$ is called $r$-{\sl divisible}. Tverberg's famous theorem \cite{tv} says the following.

\begin{theorem}\label{th:Tv} Assume $X=\{a_1,\ldots,a_n\}\subset \R^d$ where $n=(r-1)(d+1)+1$ and $r\ge 2$, $d\ge 1$ are integers. Then $X$ is $r$-divisible.
\end{theorem}

This is very simple when $d=1$. The case $r=2$ is Radon's theorem \cite{Rad}. From now one we assume that $d \ge 2$ and $r\ge 3$ are integers and $n=(r-1)(d+1)+1$.  It follows from the second part of Lemma~\ref{l:cap} (see below) that in general $X$ is not $r$-divisible when $|X|<n$.

\smallskip
Tverberg's theorem has several proofs, two by Tverberg himself \cite{tv}, \cite{tv81}, one by Tverberg and Vre\'cica~\cite{TvVrec}, by Sarkaria~\cite{Sar}, by B\'ar\'any and Onn \cite{BarOnn}, by Roudneff~\cite{Roud}, and by Zvagelski \cite{Zva}, none of them easy. The aim of this paper is to give another proof, based on linear algebra, specially structured matrices, and the theory of linear equations. It uses Tverberg's original method of ``moving the points''. For this end we have to assume that the points of $X$ are in suitably general position, namely, we assume that the coordinates of the $a_i$ (altogether $dn$ real numbers) are algebraically independent. In this case we say that $X$ is in {\sl a-general position}.

\begin{theorem}\label{th:Tva-gen} If $X\subset \R^d$ with $|X|=n$ is in a-general position and $n=(r-1)(d+1)+1$, then $X$ is $r$-divisible.
\end{theorem}

A standard limiting argument (present in Tverberg's original paper) shows that Theorem~\ref{th:Tva-gen} implies Tverberg's theorem.

\smallskip
A partition $X=\{X_1,\ldots,X_r\}$ of $X$ is called {\sl proper} if $1\le |X_p|\le d+1$ for every $p \in [r]$. Here and in what follows $[r]$ stands for the set $\{1,\ldots,r\}$. The advantage of the a-general position is shown by the following result, for a proof see for instance \cite{Bar} or \cite{BarSob}.

\begin{lemma}\label{l:cap}If $X\subset \R^d$ with $|X|=n$ is in a-general position and $\p=\{X_1,\ldots,X_r\}$ is a proper partition of $X$, then the intersection of the affine hulls of the $X_p$s is a single point $z$, that is, $z=\bigcap_{p=1}^r \aff X_p$. If, under the same conditions $|X|<n$, then $\bigcap_{p=1}^r \aff X_p=\emptyset$.
\end{lemma}

The first part of the lemma says that the following system of linear equations has a unique solution:
\begin{equation}\label{eq:aff}
z=\sum_{x\in X_p}\al(x)x \mbox{ and }1=\sum_{x\in X_p}\al(x) \mbox{ for all }p \in [r].
\end{equation}

In this form Tverberg's theorem states that there is a partition $\p$ of $X$ such that $\al(x)\ge 0$ for every $x\in X$ in the solution to  (\ref{eq:aff}).

\section{The matrix form of (\ref{eq:aff})}\label{sec:matrix}

Note that the partition $\p$ of $X$ defines a partition $J_1,\ldots,J_r$ of $[n]$ via $J_p=\{i\in [n]: a_i\in X_p\}.$ One piece of notation: for $a \in \R^d$ we denote by $(a,1)$ the vector in $\R^{d+1}$ whose first $d$ coordinates coincide with those of $a$ and the last one is $1$. We write equation (\ref{eq:aff}) in matrix form
\begin{equation}\label{matrix}
M(\al,z) =c.
\end{equation}

Here the $(n+d)\times(n+d)$ matrix $M$ is made up of blocks. The block corresponding to $X_p$ is a $(d+1)\times |J_p|$ matrix $B_p$ whose column $i \in J_p$ is $(a_i,1)$, and this block is in rows $(p-1)(d+1)+1,\ldots,p(d+1)$. The last row (of all ones) of $B_p$ is shown in Table 1.
\begin{table}[h]\label{fig:1}
\begin{center}
\begin{tabular}{|c|c|c|c|c|}
  \hline
   & & & & \\
   $B_1$ & & & & \hspace*{.1in} $-I_d$
                           \hspace*{.1in}  \\
   & & & & \\
   \hline
     1 $\dots$ 1 & & & & \\
  \hline
   & & & & \\
  &  $B_2$ & & & $-I_d$ \\
   & & & & \\
   \hline
    &  1 $\dots$ 1 & & & \\
  \hline
  & & $\ddots$ & & \\
  \hline
   & & & & \\
   & & & $B_r$ & $-I_d$ \\
   & & & & \\
  \hline
     & & &  1 $\dots$ 1  & \\
     \hline
\end{tabular}
\bigskip
\end{center}
\caption{The matrix $M$, the empty regions indicate zeros}
\end{table}
There are $r$ further blocks, each one is $-I_d$, the negative $d \times d$ identity matrix. They are in the last $d$ columns of $M$, with a row of zeroes between them. These submatrices are arranged in $M$ as shown in Table 1. All other entries of $M$ are zeroes. Column $i$ of $M$ corresponds to the vector $a_i$. Note that $M=M(\p)$ depends on $X$ and on the partition $\p=\{X_1,\ldots,X_r\}$ as well. The variables are $(\al,z) =(\al_1,\ldots,\al_n,z_1,\ldots,z_d) \in \R^{n+d}.$ The coordinates of right hand side vector $c\in \R^{n+d}$ are zero everywhere except in positions $d+1,2(d+1),\ldots,r(d+1)$ where they are equal to one.

Let $M_i$ denote the matrix obtained by replacing column $i$ of $M$ by the vector $c$. We will need the following fact.

\begin{prop}\label{base} If the partition of $X$ is proper, then $\det M\ne 0$ and $\det M_i \ne 0$ for all $i \in [n]$.
\end{prop}

{\bf Proof.} The system (\ref{matrix}) has a unique solution iff $\det M\ne 0$, which happens iff $\bigcap_{p=1}^r \aff X_p$ is a single point. So $\det M\ne 0$ is implied by Lemma~\ref{l:cap}. Next $\det M_i =0$ implies by Cramer's rule that $\al_i=\frac{\det M_i}{\det M} =0$ and so  $\bigcap_{p=1}^r \aff (X_p\setminus \{a_i\})\ne \emptyset$ which is impossible according to second half of the same lemma.\qed

\section{Moving the points}\label{sec:moving}

The method of moving the points is similar to induction. It starts with (the construction of) an initial set in $\R^d$ of size $n$ in a-general position which is $r$-divisible. The construction is simple. Take $(r-1)$ randomly rotated copies of the regular simplex inscribed in the unit sphere of $\R^d$, their vertex sets are $V_1,\ldots,V_{r-1}$ and let $v_n$ be the origin. Then $v_n \in \bigcap_1^{r-1}\conv V_p$. Set $V=\{v_n\}\cup \bigcup_1^{r-1}V_p$ and choose a point $v^*$ in the $\delta$-neighbourhood of $v$ for every $v \in V$. Set $V_p^*=\{v^*: v\in V_p\}$. If $\delta>0$ is small enough, then $v_n^*\in \bigcap_1^{r-1}\conv V_p^*$. It is clear that the vectors $v^*$ can be chosen so that $V^*$ is in a-general position. So $V^*$ is a suitable initial set. 

\smallskip
{\bf Remark.} Assume that $B=\{b_1,\ldots,b_n\}\subset \R^d$, a set in a-general position, is given in advance. It is easy to see that the initial set $V^*$ can be chosen so that the set $B \cup V^*$ is in a-general position. 

\smallskip
The crucial part of the method of moving the points is this.

\begin{theorem}\label{th:induc} If the set $\{b_1,a_1,\ldots,a_n\}\subset \R^d$ is in a-general position and the set $\{a_1,a_2,\ldots,a_n\}$ is $r$-divisible, then so is $\{b_1,a_2,\ldots,a_n\}$.
\end{theorem}

This result together with the above remark implies Theorem~\ref{th:Tva-gen}. The proof of Theorem~\ref{th:induc} starts now and will be finished in Section~\ref{sec:Xj}.

\smallskip
We write $a_1(t)=(1-t)a_1+tb_1$ and $X(t)=\{a_1(t),a_2,\ldots,a_n\}$. Let $X(t)$ be $r$-divisible with a partition $\p(t)$ and write $M(t)=M(\p,t)$ for the matrix in (\ref{eq:aff}). Typically $\p$ will not change on an interval and we use simply $M(t)$ instead of $M(\p,t)$ when there is no danger of confusion. We write $M_i(t)$ (or $M_i(\p,t)$) for the matrix where column $i$ of $M(t)$ is replaced by $c$ and all other columns are the same.

\smallskip
A value $t_0 \in [0,1)$ $X(t)$ is called {\sl regular} if $X(t)$ is $r$-divisible for all $t\in [0,t_0]$ and the solution to the system (\ref{matrix}) with $M=M(\p(t_0),t_0)$ satisfies $\al_i(t_0)>0$ for all $i.$ For instance $t_0=0$ is regular value. We fix a regular value $t_0\in [0,1)$ together with the partition $\p=\p(t_0)$, and write $M(t)$ for $M(\p,t)$.

\begin{claim}\label{cl:det} $\det M(t_0) > 0$.
\end{claim}

{\bf Proof.} Observe that $M_1(t)$ does not depend on $t$ and then $M_1(t)$ is a constant. Because of algebraic independence this constant is non-zero, say $M_1(t)=\gamma\ne 0$ and we assume that $\gamma>0$ by swapping two columns in some block if necessary.

\smallskip
The system (\ref{matrix}) with $M=M(t_0)$ has a solution so $\rank M(t_0)= \rank [M(t_0),c]$ where $[M(t_0),c]$  is the matrix $M(t_0)$ appended with $c$ as its last column. As $[M(t_0),c]$ contains $M_1(t_0)$ whose rank is $n+d$, $\rank M(t_0)=n+d$ and then $\det M(t_0)\ne 0$. Moreover, $\al_1(t_0)=\frac {\det M_1(t_0)}{\det M(t_0)}=\frac {\gamma}{\det M(t_0)}>0$ implies $\det M(t_0)>0$. \qed

\smallskip
Consider the system $M(t)(\al,z)=c$ where the only changing column is the first one. As $\det M(t)$ and $\det M_i(t)$ are continuous functions of $t$, $\al_i(t)>0$ for all $i$ and for all $t>t_0$ close enough to $t_0$. Let $t_1>t_0$ be the smallest $t>t_0$ where this fails, then either $\al_h(t_1)=0$ for some $h \in [n]$ or $\det M(t_1)=0$. Call such a $t_1$ {\sl singular}. We are done if $t_1\ge 1$, so suppose $t_1<1.$

\smallskip
We check first that $\det M(t_1)\ne 0$. As $\det M(t)$ tends to $\det M(t_1)$ when $t \to t_1^-$, $\det M(t_1)=0$ would imply that $\al_1(t)=\frac {\det M_1(t)}{\det M(t)}=\frac {\gamma}{\det M(t)}>1$ if $t<t_1$ is close enough to $t_1$ which is impossible.

\begin{claim}\label{cl:huniq} $\al_h(t_1)=0$ happens for a unique $h \in [n]$.
\end{claim}

{\bf Proof.} Since $\det M_i(t)=(1-t)\det M_i(0)+t\det M_i(1)$ is a linear function of $t$, $\det M_i(t_1)=0$ and $\det M_h(t_1)=0$ imply that
\[
t_1= \frac {\det M_i(0)}{\det M_i(0)-\det M_i(1)}=\frac {\det M_h(0)}{\det M_h(0)-\det M_h(1)}.
\]
The last equality shows $\det M_i(0)\det M_h(1)=\det M_h(0)\det M_i(1)$ contradicting algebraic independence.\qed

\smallskip
{\bf Remark.} This proof shows that the number of singular values is finite. Indeed, a singular value $t^*$ is associated with a (proper) partition $\p$ of $X(t^*)$ and an $h \in [n]$ with $\al_h(t^*)=0$ implying 
$$
t^*=\frac {\det M_h(\p,0)}{\det M_h(\p,0)-\det M_h(\p,1)},
$$
and the number of such pairs $(\p,h)$ is finite.

\smallskip
As $\det M_h(t)$ changes sign at $t=t_1$ the partition $\p$ will not work for $t>t_1$ (close to $t_1$) because there $\al_h(t)<0$. We record here that $\det M_h(t)<0$ for $t>t_1$. Assume $a_h$ is in piece $X_q$ of the partition $\p$. We will remove $a_h$ from $X_q$ and put it in some other $X_p$ to get a new partition $\p^*$ as explained next.

\section{Finding the new partition}\label{sec:Xj}

We denote by $L(t_1,\eps)$ the interval $[t_1,t_1+\eps]$ where $\eps>0$ always.
Let $\overline{M}(t)$ be the same matrix as $M(t)$ except that its column $h$ is replaced by the vector $(a_h,1,a_h,1\ldots,a_h,1)$ where $a_h,1$ is repeated $r$ times. The system $\overline{M}(t_1)(\al,z)=c$ has more than one solution: one is just the solution to $M(t_1)(\al,z)=c$ (with $\al_h=0$), the other one is $\al_h=1$ and all other $\al_i=0$ and $z=a_h$. This implies that $\det \overline{M}(t_1)=0.$

\smallskip
Let $\p_p$ be the partition you get by moving $a_h$ from $X_q$ to $X_p$ and let $\overline{M^p}(t)$ the matrix corresponding to $\p_p$. Of course $\p_q=\p$ and $\overline{M^q}(t)=M(t).$ 
Set $J=\{p\in [r]\setminus \{q\}: |X_p|\le d\}.$ Expanding $\overline{M}(t_1)$ along column $h$ gives
\begin{eqnarray*}
0=\det \overline{M}(t_1)=\sum_{p=1}^r \det \overline{M^p}(t_1)
                                     =\det M(t_1)+\sum_{p\in J} \det \overline{M^p}(t_1)
\end{eqnarray*}
because for $p\notin J$, $\overline{M^p}(t_1)$  contains $d+2$ linearly dependent columns, namely the ones in block $B_p$ (of $\p_p$) and then $\det \overline{M^p}(t_1)=0.$

\smallskip
As $\det M(t_1)>0$, $\det M(t)>0$ for all $t \in L(t_1,\eps_1)$ for some $\eps_1>0$. So there is a $j\in J$ with $ \det \overline{M^j}(t)<0$ for all $t \in L(t_1,\eps_2)$ with some $\eps_2\in (0,\eps_1]$. Fix such a $j \in J$ and set $M^*(t)=\overline{M^j}(t)$ and  $\p^*=\p_j$. This is going to be the new partition. To see this consider the system
\[
M^*(t)(\be,z)=c \mbox{ when } t \in L(t_1,\eps_2).
\]
It has a unique solution since $\rank M^*(t)=n+d$ on the interval $ L(t_1,\eps_2).$ On the same interval $\det M^*(t)<0$ and $\be_h(t)>0$ since $\det M_h^*(t)=\det M_h(t)<0$ as we recorded above. Further, for $i\in [n]$ and $i\ne h$, $\be_i(t_1)=\al_i(t_1)>0$. By continuity, $\be_i(t)$ remains positive for $t\in L(t_1,\eps_3)$ for some $\eps_3\in (0,\eps_2]$: the coefficients $\be_i(t)$ are all positive. This shows that $X(t)$ is $r$-divisible with partition $\p^*$ on $[t_1,t_2]$ where $t_2=t_1+\eps_3$, and  $t_2>t_0$ is another regular value.

\smallskip
The above proof shows that starting from $t_2$ we can find another regular value $t_4>t_2$ with a single singular value $t_3$ between them, and starting from $t_4$ another regular value $t_6$ with a singular value between them, etc. According to the remark after Claim~\ref{cl:huniq} the number of singular values is finite, and then $X(t)$ is $r$-divisible for every $t\in [0,1]$.

\bigskip
{\bf Acknowledgements.} Thanks are due to an anonymous referee for careful reading and useful advice. This piece of research was partially supported by Hungarian National Research Grants No 131529, 132696, and 133819.

\bigskip

\noindent
Imre B\'ar\'any \\
Alfr\'ed R\'enyi Institute of Mathematics,\\
13-15 Re\'altanoda Street, Budapest, 1053 Hungary\\
e-mail: {\tt imbarany@gmail.com}\\
and\\
Department of Mathematics\\
University College London\\
Gower Street, London, WC1E 6BT, UK

\end{document}